\documentclass[12pt]{article}
\usepackage{amsthm,amsmath,amssymb,latexsym,bbm}
\def\lam{\lambda}

\def\bbp{{\mathbb P}}
\def\bbe{{\mathbb E}}

\def\bbr{{\mathbb R}}

\def\call{{\cal L}}
\newtheorem{lem}{Lemma}

\def\var{\mathop{\rm Var}}
\def\lip{\mathop{\rm Lip}}

\begin{document}

\title{A Covariance Representation and an Elementary Proof of the Gaussian Concentration Inequality}
\author{Christian Houdr\'e}
\maketitle

\baselineskip=20pt
These simple notes provide a self-contained elementary proof of the 
Gaussian concentration inequality which asserts that if $X\sim N(\mu,\Sigma)$ and if
$f:\bbr^d\to\bbr$ is a Lipschitz function (with respect to the Euclidean distance), then
\begin{equation}\label{eq1}
\bbp (f(X)-\bbe f(X)\ge x)\le e^{-\frac{x^2}{2\|\Sigma^{1/2}\nabla f\|^2_\infty}},\qquad
\forall\; x>0
\end{equation}
where
\begin{equation}\label{eq2}
\|\Sigma^{1/2}\nabla f\|^2_\infty:=
\sup_{x\in\bbr^d}\|\Sigma^{1/2}\nabla f(x)\|^2=\sup_{x\in\bbr^d}\langle\nabla f(x),
\Sigma\nabla f(x)\rangle.  
\end{equation}
(Throughout $\|\cdot\|$ is the Euclidean norm, $\langle\cdot,\cdot\rangle$ is the 
Euclidean inner product, $\Sigma^{1/2}$ is a positive 
semi-definite square root of $\Sigma$, and $\nabla$ is the gradient operator).

The proof of \eqref{eq1} presented here depends on a covariance representation which
follows from interpolating characteristic functions, i.e., {\em interpolating in the Fourier domain}.  
This is formalized as:

\begin{lem}\label{lem1}
In $\bbr^d$, let $X\sim N(\mu,\Sigma)$ and let $f,g:\bbr^d\to \bbr$ be such that
$\|\nabla f\|_\infty < +\infty$ and $\|\nabla g\|_\infty <+\infty$.
Then
\begin{equation}\label{eq14}
\bbe f(X)g(X)-\bbe f(X)\bbe g(X)=\int^1_0 \bbe \langle\Sigma\nabla f(X_\alpha),\nabla g(Y_\alpha)\rangle d\alpha, 
\end{equation}
where
\begin{equation}\label{eq15}
\begin{pmatrix} X_\alpha\\ Y_\alpha\end{pmatrix}\sim N
\left(\begin{pmatrix} \mu\\ \mu\end{pmatrix}, \begin{pmatrix}
\Sigma & \alpha\Sigma\\ \alpha\Sigma & \Sigma\end{pmatrix}\right),\qquad
0\le\alpha\le 1.
\end{equation}
\end{lem}

\

\noindent
{\bf Remark~1.}  (i)  It is important to note that the first $d$-dimensional marginal of 
$\begin{pmatrix} X_\alpha\\ Y_\alpha\end{pmatrix}$ has the same law as $X$ and so does 
its second $d$-dimensional marginal.  

(ii) The representation formula \eqref{eq14} is well known (see the Bibliographical Comments 
at the end of these pages).  An elementary proof of it, 
based on characteristic functions, will be shown below.  
It can also be obtained as a direct consequence of It${\hat {\rm o}}$'s formula 
or, in a more elementary way, using the so-called Plackett's identities.
For another approach, let $\gamma_d$ denote the canonical Gaussian 
measure on $\bbr^d$ with density $\gamma_d(dx)/dx = (2\pi)^{-d/2}e^{-\|x\|^2/2}$.  Let 
$(P_t)_{t\ge 0}$ be the Ornstein-Uhlenbeck semigroup 
associated to $\gamma_d$ and acting on the space $L^1(\gamma_d)$ via 
$$(P_tg)(x)=\int_{\bbr^d}g(e^{-t}x + \sqrt{1-e^{-2t}}y)\gamma_d(dy),       x\in \bbr^d.$$ 
Clearly, for $g$ smooth, $\nabla P_tg = e^{-t} P_t\nabla g$, and an equivalent version 
of \eqref{eq14} (for $\Sigma = I_d$, the $d\times d$ identity matrix) is given by:   

\begin{equation}\label{eqOU}
\bbe fg-\bbe f\bbe g=\int^{+\infty}_0 \bbe \langle\nabla f,\nabla P_tg \rangle dt 
= \int^{+\infty}_0 e^{-t}\bbe \langle\nabla f,P_t\nabla g \rangle dt , 
\end{equation}
since if $Y$ is an independent copy of $X$, $\begin{pmatrix} X\\ e^{-t}X+ \sqrt{1-e^{-2t}}Y\end{pmatrix}$ has the same 
law as $\begin{pmatrix} X_{e^{-t}}\\ Y_{e^{-t}}\end{pmatrix}$.  
The covariance representation \eqref{eq14} which, in contrast to 
\eqref{eqOU}, hides the presence of the Ornstein-Uhlenbeck semigroup 
seems to be simpler and more amenable to our purposes, avoiding the study of actions of the semigroup on classes of functions.

\

\noindent
{\bf Proof of the concentration inequality \eqref{eq1}}

Let us now show how a simple application of \eqref{eq14} gives \eqref{eq1}.
Again, note that $X_\alpha\buildrel\call\over =X$ and
$Y_\alpha\buildrel\call\over=  Y$.  Then, taking $g=e^{tf}$, with $f$ bounded and 
such that
$\|\nabla f\|_\infty< +\infty$, with moreover $\bbe f(X)=0$, we have
\begin{align}\label{eq3}
\bbe f(X)e^{tf(X)}&= t\int^1_0\bbe \langle\nabla f(X_\alpha),
\Sigma\nabla f(Y_\alpha)e^{tf(Y_\alpha)}\rangle d\alpha\nonumber\\
&= t\int^1_0 \bbe \langle\Sigma^{1/2}\nabla f(X_\alpha),\Sigma^{1/2}\nabla 
f(Y_\alpha)\rangle e^{tf(Y_\alpha)}d\alpha\nonumber\\
&\le t\int^1_0 \bbe \|\Sigma^{1/2}\nabla f(X_\alpha)\|\| \Sigma^{1/2}\nabla
f(Y_\alpha)\|e^{tf(Y_\alpha)}d_\alpha\nonumber\\
&\le t\| \Sigma^{1/2}\nabla f\|^2_\infty\int^1_0 \bbe e^{tf(Y_\alpha)}d\alpha\nonumber\\
&=t\| \Sigma^{1/2}\nabla f\|^2_\infty \bbe e^{tf(X)}.
\end{align}
Setting $h(t)=\bbe e^{tf}$, the inequality \eqref{eq3} reads as
\begin{equation}\label{eq4}
\forall\; t>0,\qquad h'(t)\le t\|\Sigma^{1/2}\nabla f\|^2_\infty h(t),
\end{equation}
for any $f$ bounded such that $\bbe f=0$ and $\|\nabla f\|_\infty < +\infty$.
Thus,
\begin{equation}\label{eq5}
\bbe e^{t(f-\bbe f)}\le e^{{t^2\|\Sigma^{1/2}\nabla f\|^2_\infty}/{2}}
\end{equation}
for any $f$ bounded such that $\|\nabla f\|_\infty < +\infty$.  The 
Bienaym\'e-Chebyshev inequality gives, 
\begin{equation}\label{eq6}
\bbp (f-\bbe f\ge x)\le\inf_{t>0} e^{-tx}\bbe e^{t(f-\bbe f)} =
e^{-{x^2}/{2\|\Sigma^{1/2}\nabla f\|^2_\infty}},
\end{equation}
for all $x>0$.  
Finally, to remove the boundedness assumption, approximate $f$ by
\begin{equation}\label{eq7}
f_n=f\mathbbm 1_{|f|< n}+|n|\mathbbm 1_{|f| \ge n}.
\end{equation}
Then $|f_n(x)-f_n(y)|\le |f(x)-f(y)|$, $\|\Sigma^{1/2}\nabla f_n\|_\infty\le 
\|\Sigma^{1/2}\nabla f\|_\infty$,
while $\bbe f_n\to \bbe f$.  \hfill$\Box$

\

\noindent
{\bf Remark~2.}  (i) In the inequality \eqref{eq1}, we then are 
often left with the problem of estimating $\| \Sigma^{1/2}\nabla f\|^2_\infty$.
A generic estimate is obtained via the weak second moment of $X$, i.e., via 
\begin{equation}\label{eq8}
\|\Sigma^{1/2}\nabla f\|^2_\infty\le \lambda^*\|\nabla f\|^2_\infty,
\end{equation}
where $\lam^*$ is the maximal eigenvalue of $\Sigma$:
\begin{equation}\label{eq9}
\lam^*=\sup_{\|x\|=1}\langle \Sigma x,x\rangle=\sup_{\|x\|=1}\bbe|\langle X,x\rangle|^2.  
\end{equation}

(ii) In some instances, e.g., for some particular classes of functions, the weak
moment of order two of $X$ can be replaced by the strong second moment, i.e., by
$\sigma^{*2}=\max_{i=1,\dots, d}\sigma_{i,i} = \max_{i=1,\dots, d}\var X_i$.
A case at hand is $f(x)=\max_{1\le i\le d}x_i$.  
Indeed, then, 
\begin{equation}\label{eq10}
\langle \Sigma^{1/2}\nabla f(x),\Sigma^{1/2}\nabla f(x)\rangle=
\sum^d_{i=1}\sum^d_{j=1}\sigma_{i,j} \frac{\partial f (x)}{\partial x_i}
\frac{\partial f(x)}{\partial x_j}
=\sum^d_{i=1}\sigma_{i,i}\left(\frac{\partial f(x)}{\partial x_i}\right)^2,
\end{equation}
since for $i\ne j$, ${\partial f}/{\partial x_i}$ and ${\partial f}/{\partial x_j}$
live on disjoint sets.  Now \eqref{eq10} is itself upper bounded by
\begin{equation}\label{eq11}
\max_{1\le i\le d}\sigma_{i,i} \sum^d_{i=1}\left(\frac{\partial f(x)}{\partial x_i}\right)^2
\le \max_{1\le i\le d}\sigma_{i,i} = \sigma^{*2}, 
\end{equation}
since $\sum^d_{i=1}\left({\partial f}/{\partial x_i}\right)^2\le 1$.  
There are a few ways to obtain this last fact: For example,
if $e=(1,\dots, 1)\in \bbr^d$, and for $x\in\bbr^d$, $t\in \bbr$,
$\max(x+te)=(\max x)+t$.  Hence by the chain rule,
$1=d\max(x+te)/dt=\sum^d_{i=1}{\partial \max x}/{\partial x_i}$.

Another way to show \eqref{eq10} is to note that
\begin{equation}\label{eq12}
|\max (x_1,\dots, x_d)-\max (y_1,\dots, y_d)|\le 
\max_{1\le i\le d} |x_i-y_i|\le\sqrt{\sum^d_{i=1}(x_i-y_i)^2}, 
\end{equation}
and then to appeal to the fact that if $f\in\lip (1)$, with respect to the Euclidean
norm, then (Rademacher Theorem) $\sum^d_{i=1}\left(\frac{\partial f}{\partial x_i}\right)^2\le 1$ .

(iii) With a bit more effort, it is possible to get an improved version of the 
Gaussian concentration inequality, i.e.,
\begin{align}
\bbp (f(X)-\bbe f(X)\ge x) 
&\le \frac{\bbe(f-\bbe f)^+}{x}e^{-x^2/2\|\Sigma^{1/2}\nabla f\|^2_\infty} \label{eq13}\\
&\le \frac{\sqrt \pi}{2\sqrt{2}x}\|\Sigma^{1/2}\nabla f\|_\infty
e^{-x^2/2\|\Sigma^{1/2}\nabla f\|^2_\infty},\label{eq23}
\end{align}
for all $x > 0$.  

(iv) The term $\|\Sigma^{1/2}\nabla f\|^2_\infty$, appears in \eqref{eq1} 
because of the use of the
Cauchy-Schwarz inequality. Different types of upper bounds occur by applying, instead, 
Young's inequality, potentially leading to improve concentration results for special 
classes of functions.

\

\noindent
{\bf Proof of Lemma~1} The idea of the proof is to interpolate between the independent
case, $\alpha=0$ and the "fully dependent identical" case $\alpha=1$.  Here is what we mean:

In $\bbr^{2d}$, we form the vector $\begin{pmatrix}X_\alpha\\ Y_\alpha\end{pmatrix}$,
$0\le \alpha\le 1$, which is Gaussian with mean vector
$\alpha\begin{pmatrix} \mu\\ \mu\end{pmatrix} + (1-\alpha)
\begin{pmatrix}\mu\\ \mu\end{pmatrix}=\begin{pmatrix}\mu\\ \mu\end{pmatrix}$
and covariance matrix $\alpha\begin{pmatrix} \Sigma & \Sigma\\ \Sigma & \Sigma\end{pmatrix}
+(1-\alpha)\begin{pmatrix} \Sigma & 0\\ 0 & \Sigma\end{pmatrix} = 
\begin{pmatrix}\Sigma & \alpha\Sigma\\ \alpha\Sigma & \Sigma\end{pmatrix}$.
So for $\alpha=0$ the vector $\begin{pmatrix}X_0\\ Y_0\end{pmatrix}$ is normal with
$X_0$ and $Y_0$ independent and $X_0\buildrel \call\over = Y_0\buildrel \call\over = X$,
while for $\alpha=1$, $\begin{pmatrix} X_1\\ Y_1\end{pmatrix}$ is a Gaussian vector on
$\bbr^{2d}$ with $X_1=Y_1$ and $X_1\buildrel \call\over = Y_1\buildrel \call\over = X$.
Hence if $\bbe _\alpha$ denotes the expectation with respect to the Gaussian
measure corresponding to $\begin{pmatrix} X_\alpha\\ Y_\alpha\end{pmatrix}$, we have
\begin{align}\label{eq16}
\bbe (f(X)g(X))&= \bbe _1(f\otimes g)\nonumber\\
\bbe f(X)\bbe g(X)&= \bbe _0(f\otimes g),
\end{align}
where 
\begin{align*}f\otimes g: \bbr^d\times\bbr^d&\longrightarrow\bbr\\
(x,y)&\longrightarrow (f\otimes g)(xy)=f(x)g(y).\end{align*}
Now the characteristic function of $\begin{pmatrix}X_0\\ Y_0\end{pmatrix}$ is given by
\begin{equation}\label{eq17}
\varphi_0(t,s):=\varphi_{(X_0,Y_0)}(t,s)=\varphi_X(t)\varphi_X (s),
\end{equation}
where $\varphi_X$ is the characteristic function of $X$, i.e., 
$\varphi_X(t)=e^{i\langle\mu, t\rangle -\frac12 \langle\Sigma t,t\rangle}$.  
Now, the characteristic
function of $\begin{pmatrix}X_1\\ Y_1\end{pmatrix}$ is given by
\begin{equation}\label{eq18}
\varphi_1(t,s):=\varphi_{(X_1,Y_1)}(t,s)=\varphi_X(t+s).
\end{equation}
Moreover, 
\begin{align}\label{eq19}
\varphi_{(X_\alpha,Y_\alpha)}(t,s)&=\bbe e^{i\langle (t,s),(X_\alpha,Y_\alpha)\rangle}
= \left[ \varphi_{(X_1,Y_1)} (t,s)\right]^\alpha
\left[ \varphi_{(X_0,Y_0)} (t,s)\right]^{1-\alpha}\nonumber\\
&=\varphi^\alpha_X(t+s)\varphi^{1-\alpha}_X(t)\varphi^{1-\alpha}_X(s).  
\end{align}
So using the notation $\varphi_\alpha (t,s)$, $0\le \alpha\le 1$, to prove
\eqref{eq14} note that 
\begin{equation}\label{eq20}
\varphi_1(t,s)-\varphi_0 (t,s)=\int^1_0 \frac d{d\alpha}\varphi_\alpha (t,s)d\alpha,
\end{equation}
and that
\begin{equation}\label{eq21}
\frac d{d\alpha}\varphi_\alpha(t,s)= -\langle\Sigma t,s\rangle \varphi_\alpha(t,s) = \bbe \langle\Sigma\nabla 
f(X_\alpha),\nabla g(Y_\alpha)\rangle,  
\end{equation}
for the functions
$f(x)=e^{i\langle t,x\rangle}$ and $g(x)=e^{i\langle s,x\rangle}$.  
Plugging \eqref{eq21} into \eqref{eq20} just gives \eqref{eq14} for such $f$ and $g$.   
Then by bilinearity, the equality \eqref{eq14} extends to trigonometric
polynomial $\sum^N_{i=1}a_je^{i\langle t_j,x\rangle}$ and 
$g(x)=\sum^N_{j=1} b_je^{i\langle s_j,x\rangle}$  from there ``by density" to the 
class of functions described in the lemma. \hfill$\Box$

\

{\bf Bibliographical Comments:}  As already indicated, 
the representation formula \eqref{eq14} is well known.  We refer 
the reader to:  Houdr\'e, C.; P\'erez-Abreu, V.; Surgailis, D.,  
Interpolation, correlation identities, and inequalities for infinitely divisible variables.  
{\it J. Fourier Anal. Appl.}  {\bf 4}  (1998),  651--668, 
for precise references and connections as well as for 
a complete proof of Lemma~1 along the lines given above.    
The proof of \eqref{eq1} presented here is in: Bobkov, S.~G.; G\"otze, F.; Houdr\'e, C.;  
On Gaussian and Bernoulli covariance representations.  
{\it Bernoulli} {\bf 7}  (2001),  439--451., where the reader will also find a 
proof of \eqref{eq13} and \eqref{eq23} (with wrong absolute constants!).  
An extension of \eqref{eq14} 
to infinitely divisible (ID) laws as well as further applications 
can be found in:  Houdr\'e, C.; Comparison and deviation 
from a representation formula.  {\it Stochastic processes and related topics,}  207--218, Trends Math., 
Birkh\"auser Boston, Boston, MA, 1998.   In turn this extension allows for various concentration inequalities in the ID framework.  

\end{document}